\documentclass[11pt]{article}

\usepackage{amsmath,amssymb,amsthm}
\usepackage{hyperref}

\newtheorem{theorem}{Theorem}[section]
\newtheorem{lemma}[theorem]{Lemma}

\newtheorem{remark}[theorem]{Remark}

\newtheorem{corollary}[theorem]{Corollary}

\title{Global existence for the 1D compressible Euler equations with time-decaying damping and spatial perturbations}
\author{
Yuusuke Sugiyama\\
Department of Mathematics, Tokyo University of Science\\
Kagurazaka 1-3, Shinjuku-ku, Tokyo 162-8601, Japan\\
\texttt{sugiyama.y@rs.tus.ac.jp}
}
\date{}

\begin{document}
\maketitle %

\begin{abstract}
We study the Cauchy problem for the one-dimensional
compressible Euler equations with damping coefficient
$b(t)+a(x)$, where
$b(t)=\lambda(1+t)^{-\mu}$ with
$0<\mu<1$, $\lambda>0$, or $\mu=1$, $\lambda>2$,
and $a\in C_b^1(\mathbb R)\cap L^1(\mathbb R)$ is
nonnegative.

We prove global existence and uniqueness of classical
solutions for sufficiently small $C_b^1$ perturbations
of a constant state, without assuming spatial decay or
convergence of the initial data at infinity.
No smallness condition is imposed on the spatially
dependent damping coefficient, even though its
contribution persists for all time.

We also obtain decay estimates of order
$\varepsilon(1+t)^{-\mu}$ for the derivatives
$u_t$, $v_t$, and $v_x$, together with a pointwise
estimate for $u_x$ involving the spatial coefficient
$a(x)$.

The main ingredient is a time-uniform $L^\infty$
estimate for an inhomogeneous coupled transport system
with a non-small spatial coefficient.
Combined with suitable transformations of the
Riemann-invariant derivatives, this estimate allows
us to control the nonlinear terms and establish
global $C^1$ bounds.
\end{abstract}

\section{Problem and main result}

\subsection{Cauchy problem}

We consider the Cauchy problem for the one-dimensional
compressible Euler equations with time-dependent damping
and a spatial perturbation:
\begin{equation}
\begin{cases}
u_t-v_x=0,
\\[1mm]
v_t+p(u)_x=-\bigl(b(t)+a(x)\bigr)v,
\\[1mm]
(u(0,x),v(0,x))
=(1+\varepsilon\varphi(x),\varepsilon\psi(x)),
\end{cases}
\qquad (t,x)\in[0,\infty)\times\mathbb R.
\label{eq:Euler}
\end{equation}

Here $u=u(t,x)>0$ denotes the specific volume,
$v=v(t,x)$ denotes the fluid velocity, and
$\varepsilon>0$ is a small parameter.

We assume that the pressure is given by
\begin{equation}
p(u)=\frac{u^{-\gamma}}{\gamma},
\qquad \gamma>1.
\label{eq:pressure}
\end{equation}

The sound speed is defined by
\begin{equation}
c(u)=\sqrt{-p'(u)}
=u^{-\frac{\gamma+1}{2}}.
\label{eq:sound-speed}
\end{equation}

We assume that
\begin{equation}
\varphi,\psi\in C_b^1(\mathbb R).
\label{eq:initial-data}
\end{equation}

Here $C_b^1(\mathbb R)$ denotes the space of bounded
continuously differentiable functions whose first
derivatives are also bounded.

\subsection{Assumptions on the damping coefficients}

We assume that the time-dependent damping coefficient
$b(t)$ is given by
\begin{equation}
b(t)=\frac{\lambda}{(1+t)^\mu},
\qquad t\geq0,
\label{eq:time-damping}
\end{equation}
where the parameters $\lambda$ and $\mu$ satisfy
\begin{equation}
\begin{cases}
\lambda>0,\quad 0<\mu<1,
\\
\text{or}
\\
\lambda>2,\quad \mu=1.
\end{cases}
\label{eq:parameters}
\end{equation}

For the spatial perturbation $a(x)$, we impose
\begin{equation}
a\in C_b^1(\mathbb R),
\qquad
a(x)\geq0
\quad\text{for all }x\in\mathbb R,
\label{eq:a-regularity}
\end{equation}
and
\begin{equation}
\int_{\mathbb R}a(x)\,dx<\infty.
\label{eq:a-integrability}
\end{equation}

We emphasize that no smallness assumption is imposed
on $a(x)$ or its $L^1$-norm.

\subsection{Main result}

Our main result is the global existence of small
classical solutions to \eqref{eq:Euler}.

\begin{theorem}[Global existence and decay]
\label{thm:global}

Let $\gamma>1$, and suppose that
\eqref{eq:initial-data}--\eqref{eq:a-integrability}
are satisfied.

Then there exists a constant $\varepsilon_0>0$
such that, for every $0<\varepsilon\leq\varepsilon_0$,
the Cauchy problem \eqref{eq:Euler}
admits a unique global classical solution
\[
(u,v)\in C_b^1([0,\infty)\times\mathbb R)
\]
satisfying $u(t,x)>0$ for all $t\geq0$ and
$x\in\mathbb R$.

Moreover, the solution satisfies
\begin{equation}
\|u(t,\cdot)-1\|_{L^\infty}
+
\|v(t,\cdot)\|_{L^\infty}
\leq C\varepsilon
\label{eq:zeroth-estimate}
\end{equation}
for all $t\geq0$.

Furthermore, we have
\begin{equation}
\|u_t(t,\cdot)\|_{L^\infty}
+
\|v_t(t,\cdot)\|_{L^\infty}
+
\|v_x(t,\cdot)\|_{L^\infty}
\leq
C\varepsilon(1+t)^{-\mu}
\label{eq:time-derivative-estimate}
\end{equation}
and
\begin{equation}
|u_x(t,x)|
\leq
C\varepsilon
\left(
(1+t)^{-\mu}+a(x)
\right)
\label{eq:space-derivative-estimate}
\end{equation}
for all $(t,x)\in[0,\infty)\times\mathbb R$.

The constants $\varepsilon_0$ and $C$ may depend
on $\gamma$, $\lambda$, $\mu$, the initial profiles
$\varphi,\psi$, and the spatial coefficient $a$,
but are independent of $t$ and $\varepsilon$.
\end{theorem}

\begin{remark}
The initial profiles $\varphi$ and $\psi$
are not assumed to have compact support or
to converge as $|x|\to\infty$.

Furthermore, the spatial coefficient $a(x)$
is not required to be small.
For example, the assumptions allow
\[
a(x)=\frac{a_0}{(1+x^2)^{\nu/2}},
\qquad
a_0\geq0,\quad \nu>1,
\]
where $a_0$ may be arbitrarily large.
\end{remark}

\begin{remark}
The local existence and uniqueness of classical solutions to \eqref{eq:Euler} follow from the standard theory of quasilinear hyperbolic systems (see, for example, \cite{Friedrichs1948}). More precisely, if the initial data satisfy
\[
u_0,v_0\in C_b^1(\mathbb R),
\qquad
\inf_{x\in\mathbb R}u_0(x)>0,
\]
then there exists $T>0$ such that \eqref{eq:Euler} admits a unique classical solution $(u,v)\in C_b^1([0,T]\times\mathbb R)$ with $u(t,x)>0$. Neither smallness nor spatial decay of the initial data is required for local existence.
\end{remark}

\subsection{Known results}
\label{subsec:known-results}

We first recall some results for the one-dimensional
compressible Euler equations without damping and
with constant damping.
In the absence of damping, smooth solutions can
develop singularities in finite time due to nonlinear
compression; see, for example,
\cite{Lax1964,ChenPanZhu2017}.
For a positive constant damping coefficient,
Hsiao and Liu
\cite{HsiaoLiu1992,HsiaoLiu1993}
established global existence of small smooth solutions
and their asymptotic behavior governed by a nonlinear
diffusion equation.
Further developments on global solutions and
diffusive asymptotics include
\cite{HsiaoSerre1996,Nishihara1997,
MarcatiNishihara2003,Mei2010}.

We next review the case of damping depending only
on time,
\begin{equation}
b(t)=\frac{\lambda}{(1+t)^\mu},
\qquad \lambda>0,\quad \mu\geq0.
\label{eq:known-time-damping}
\end{equation}
Pan \cite{Pan2016Global,Pan2016Blowup,Pan2021}
studied the global existence and finite-time blow-up
of small-data solutions to the compressible Euler
equations with time-dependent damping.
In the notation of \eqref{eq:known-time-damping},
global classical solutions exist for sufficiently
small initial data when
\begin{equation}
0\leq\mu<1,\quad \lambda>0,
\qquad\text{or}\qquad
\mu=1,\quad\lambda>2.
\label{eq:known-global-range}
\end{equation}
On the other hand, finite-time singularities may
occur for suitable small initial data when
$\mu>1$, or when $\mu=1$ and $0<\lambda\leq2$.

The author \cite{Sugiyama2018} obtained sharp upper
and lower lifespan estimates in the
singularity-forming regimes.
For suitable compressive initial profiles of size
$\varepsilon$, these can be summarized as
\begin{equation}
T^*(\varepsilon)\asymp
\begin{cases}
\varepsilon^{-1}, & \mu>1,\\
\varepsilon^{-2/(2-\lambda)},
& \mu=1,\ 0<\lambda<2,
\end{cases}
\label{eq:known-lifespan}
\end{equation}
while
$\log T^*(\varepsilon)\asymp\varepsilon^{-1}$
when $\mu=1$ and $\lambda=2$. Here $T^*(\varepsilon)$ denotes the maximal existence
time of the classical solution to the Cauchy problem
\eqref{eq:Euler}.

Subsequently, the author \cite{Sugiyama2020}
proved global existence in the range
\eqref{eq:known-global-range} for small
$C_b^1(\mathbb R)$ perturbations, without assuming
the existence of limits of the initial data as
$|x|\to\infty$.
Other results concerning large initial data
and the critical asymptotic behavior can be found in
\cite{ChenLiLiMeiZhang2020,GengLinMei2020,Chen2024}.

The distinction between effective and non-effective
time-dependent damping is also reflected in the
theory of linear damped wave equations.
Wirth \cite{Wirth2004,Wirth2006,Wirth2007}
investigated solution representations and long-time
behavior for such equations, including the critical
scale of the damping.

When the damping coefficient depends on the spatial
variable, new difficulties arise in controlling
the Riemann invariants along characteristic curves.
Related characteristic methods for spatially
inhomogeneous Euler systems were developed by
Chen, Young and Zhang \cite{ChenYoungZhang2013}
in the non-isentropic setting.

For purely space-dependent damping,
Geng, Lai, Yuen and Zhou
\cite{GengLaiYuenZhou2023}
established finite-time blow-up under spatial
integrability assumptions on the damping and
compact-support assumptions on the initial
perturbations.

Chen et al. \cite{ChenLiLiMeiZhang2020}
established a maximum principle for the Riemann invariants
under the assumption that the damping coefficient
$\alpha(t,x)\in C_b^1$ is positive.
They also proved global existence for damping coefficients
satisfying
\[
\begin{aligned}
C_1(1+t)^{-\nu}
&\leq \alpha(t,x)
\leq C_2(1+t)^{-\nu},\\
|\alpha_t(t,x)|+|\alpha_x(t,x)|
&\leq C(1+t)^{-2\nu},
\end{aligned}
\]
where $0<\nu<1$, under suitable assumptions on the
initial data, including their convergence at spatial infinity.

More recently, the author \cite{Sugiyama2025} considered damping coefficients $d(t,x)$ satisfying
\[
|d(t,x)-\mu_0|+|d_t(t,x)|+|d_x(t,x)|
\leq a_1(t)+a_2(x),
\]
where $\mu_0\geq0$, $a_1\in L^1([0,\infty))$ and $a_2\in L^1(\mathbb R)$ are nonnegative functions, together with some additional assumptions. The author proved global existence for sufficiently small initial data when $\mu_0>0$ and finite-time blow-up for suitable small initial data when $\mu_0=0$.

Cai, Liu, Mei and Wang
\cite{CaiLiuMeiWang2025,CaiLiuMeiWang2026}
studied the compressible Euler equations with
spatiotemporal damping of the form
\[
\frac{\alpha_1(x)}{(1+t)^\nu}=\frac{\alpha_\infty}{(1+t)^\nu}+\frac{\alpha_1(x)-\alpha_\infty}{(1+t)^\nu},
\]
where $\alpha_1(x)$ converges to a positive constant
$\alpha_\infty$ as $|x|\to\infty$.

In \cite{CaiLiuMeiWang2025}, they considered
$0\leq\nu<1$ under suitable weighted integrability
and regularity assumptions on $\alpha_1(x)$.
In \cite{CaiLiuMeiWang2026}, they studied
the critical case $\nu=1$, assuming
$\alpha_\infty>2$ and
\[
\alpha_1(x)\geq\alpha_\infty,
\qquad
\alpha_1-\alpha_\infty\in L^1(\mathbb R)\cap L^2(\mathbb R),
\qquad
\alpha_1'\in L^2(\mathbb R).
\]
They established global existence and convergence
to asymptotic profiles under suitable assumptions
on the initial data, including their convergence
to prescribed states as $x\to\pm\infty$.

The damping coefficient considered in the present
paper is of a different nature from those treated
in the preceding works.
We consider the additive coefficient
$b(t)+a(x)$, where $b(t)$ tends to zero but
is not integrable over $[0,\infty)$, whereas
$a(x)$ is independent of time and integrable
in space.

In particular, the perturbation framework of
\cite{Sugiyama2020}, which is based on a constant damping
coefficient and integrable perturbations in time
and space, does not directly apply.
Indeed, $b(t)$ cannot be treated as an
integrable perturbation of zero, and the
coefficient has no strictly positive constant
part.

When $\nu=0$, the damping coefficient in
\cite{CaiLiuMeiWang2025} is a positive constant
plus a spatial perturbation, which is also covered
by \cite{Sugiyama2025} under suitable integrability
assumptions. When $\nu>0$, both the constant and
spatially dependent parts decay in time at the
same rate. In contrast, we consider $b(t)+a(x)$,
where $b(t)$ tends to zero but is not integrable
in time, while $a(x)$ is independent of time
and integrable in space.

Although $a(x)\geq 0$ contributes additional
damping, its spatial dependence creates terms
involving $a'(x)v$ when the equations for the
Riemann invariants are differentiated.
Consequently, uniform estimates for the
Riemann invariants alone do not directly
provide the derivative bounds needed for
global classical solutions.
This difficulty is particularly relevant
when the spatial coefficient is not assumed
to be small.

The present paper resolves this difficulty
by combining transformations of the
Riemann-invariant derivatives with a
time-uniform estimate for an inhomogeneous
linear transport system.
As a consequence, global existence holds
throughout the effective time-dependent
damping range considered here, without
smallness of $\|a\|_{L^1}$ and without
spatial decay assumptions on the initial data.

\subsection{Idea of the proof}

The key idea of the proof is to introduce suitable transformations of the derivatives of the Riemann invariants and to establish uniform estimates for a linear transport system.

Differentiating the equations for the Riemann invariants produces terms involving $a'(x)v$ and cross terms involving $b(t)$. To eliminate these terms, we first introduce $R$ and $S$ by subtracting suitable correction terms from $r_x$ and $s_x$, and then define
\[
\widetilde R=R+\frac{b(t)}{2}(u-1),
\qquad
\widetilde S=S+\frac{b(t)}{2}(u-1).
\]
These transformations eliminate the terms involving $a'(x)v$ and the cross terms involving $b(t)$, such as $b(t)S/2$ in the equation for $R$. As a result, the equations for $\widetilde R$ and $\widetilde S$ can be written as a linear transport system with suitable source terms.

Chen et al. \cite{ChenLiLiMeiZhang2020} established
a maximum principle for the Riemann invariants
by introducing a suitable energy functional.
A related comparison principle was also used in
\cite{Sugiyama2025} to obtain uniform estimates
for the Riemann invariants themselves.
In the present paper, we further establish
uniform $L^\infty$ estimates for an inhomogeneous
linear transport system by constructing
time-independent supersolutions and subsolutions.
These estimates allow us to control the source
terms involving $a(x)$ without assuming
the smallness of $\|a\|_{L^1}$.
Applying these estimates to the transformed
derivatives $\widetilde R$ and $\widetilde S$,
we obtain the derivative bounds needed
for global existence.

\section{Preliminaries}

\subsection{Riemann invariants}

We first introduce the Riemann invariants associated
with the system \eqref{eq:Euler}.

Recall that
\[
c(u)=\sqrt{-p'(u)}
=u^{-\frac{\gamma+1}{2}}.
\]
We define
\[
\eta(u)
=\int_u^\infty c(\xi)\,d\xi
=\frac{2}{\gamma-1}u^{-\frac{\gamma-1}{2}}
\]
and
\[
r=v-\eta(u)+\frac{2}{\gamma-1},
\qquad
s=v+\eta(u)-\frac{2}{\gamma-1}.
\]

We denote the characteristic differential operators by
\[
\partial_-=\partial_t-c(u)\partial_x,
\qquad
\partial_+=\partial_t+c(u)\partial_x.
\]

Then the system \eqref{eq:Euler} is equivalent to
\begin{equation}
\begin{cases}
\displaystyle
\partial_-r=-\frac{b(t)+a(x)}{2}(r+s),
\\[2mm]
\displaystyle
\partial_+s=-\frac{b(t)+a(x)}{2}(r+s).
\end{cases}
\label{eq:riemann-system}
\end{equation}

Since $r+s=2v$, we also have
\begin{equation}
r_x=u_t+c(u)u_x,
\qquad
s_x=u_t-c(u)u_x.
\label{eq:riemann-derivatives}
\end{equation}

For later use, we denote by
$x_\pm(\tau;t,x)$ the characteristic curves satisfying
\[
\begin{cases}
\displaystyle
\frac{d}{d\tau}x_\pm(\tau;t,x)
=\pm c(u(\tau,x_\pm(\tau;t,x))),
\\[2mm]
x_\pm(t;t,x)=x.
\end{cases}
\]

\subsection{Estimates for a linear transport system}

In this subsection, we establish uniform estimates
for a linear transport system.
These estimates will play an important role in the
proof of the global existence theorem.

Let $T>0$ and suppose that
$c=c(t,x)$ is a given sufficiently regular function
satisfying
\begin{equation}
0<c_*\leq c(t,x)\leq c^*
\label{eq:linear-speed}
\end{equation}
for all $(t,x)\in[0,T]\times\mathbb R$.

We consider the linear system
\[
\begin{cases}
\displaystyle
\partial_-y=\frac{a(x)}{2}(z-y)+f,
\\[2mm]
\displaystyle
\partial_+z=\frac{a(x)}{2}(y-z)+g,
\end{cases}
\]
where
\[
\partial_-=\partial_t-c(t,x)\partial_x,
\qquad
\partial_+=\partial_t+c(t,x)\partial_x.
\]

For a vector-valued function $F=(f,g)$, we write
\[
\|F\|_\infty
=
\max\left\{
\|f\|_{L^\infty(\mathbb R)},
\|g\|_{L^\infty(\mathbb R)}
\right\}.
\]

We first establish a comparison principle for
a linear transport system.
Let $d\in C_b([0,T]\times\mathbb R)$ satisfy
$d(t,x)\geq0$.
We consider bounded continuous functions
$w_1,w_2$ satisfying the following integral
inequalities along the characteristic curves:
\begin{equation}
\begin{aligned}
w_1(t,x)
\leq{}&
e^{-\frac12\int_s^t d(\xi,x_-(\xi))\,d\xi}
w_1(s,x_-(s))
\\
&+\frac12\int_s^t
e^{-\frac12\int_\tau^t d(\xi,x_-(\xi))\,d\xi}
d(\tau,x_-(\tau))
w_2(\tau,x_-(\tau))\,d\tau,
\\
w_2(t,x)
\leq{}&
e^{-\frac12\int_s^t d(\xi,x_+(\xi))\,d\xi}
w_2(s,x_+(s))
\\
&+\frac12\int_s^t
e^{-\frac12\int_\tau^t d(\xi,x_+(\xi))\,d\xi}
d(\tau,x_+(\tau))
w_1(\tau,x_+(\tau))\,d\tau,
\end{aligned}
\label{eq:comparison-system}
\end{equation}
for all $0\leq s\leq t\leq T$ and $x\in\mathbb R$,
where $x_\pm(\tau)=x_\pm(\tau;t,x)$.

For continuously differentiable functions,
\eqref{eq:comparison-system} is equivalent to
\[
\begin{cases}
\displaystyle
\partial_-w_1
\leq\frac{d(t,x)}2(w_2-w_1),
\\[2mm]
\displaystyle
\partial_+w_2
\leq\frac{d(t,x)}2(w_1-w_2).
\end{cases}
\]
Thus, \eqref{eq:comparison-system} describes
subsolutions of the corresponding homogeneous
linear transport system.
The integral formulation allows us to work
with continuous functions.

\begin{lemma}[Comparison principle]
\label{lem:comparison}

Suppose that bounded continuous functions
$w_1,w_2$ satisfy \eqref{eq:comparison-system}.

If
\[
w_1(s,x)\leq0,\qquad w_2(s,x)\leq0
\]
for all $x\in\mathbb R$, then
\[
w_1(t,x)\leq0,\qquad w_2(t,x)\leq0
\]
for all $(t,x)\in[s,T]\times\mathbb R$.

In particular, if both inequalities in
\eqref{eq:comparison-system} hold with equality,
then
\begin{equation}
\begin{aligned}
\max\bigl\{
&\|w_1(t)\|_{L^\infty},
\|w_2(t)\|_{L^\infty}
\bigr\}
\\
&\leq
\max\bigl\{
\|w_1(s)\|_{L^\infty},
\|w_2(s)\|_{L^\infty}
\bigr\}
\end{aligned}
\label{eq:comparison-maximum}
\end{equation}
for all $0\leq s\leq t\leq T$.

\end{lemma}

\begin{proof}

Since $w_1(s,\cdot)\leq0$ and $d\geq0$,
the first inequality in
\eqref{eq:comparison-system} gives
\[
w_1(t,x)
\leq
\frac{\|d\|_{L^\infty}}2
\int_s^t
\|(w_2(\tau))_+\|_{L^\infty}\,d\tau.
\]
Similarly,
\[
w_2(t,x)
\leq
\frac{\|d\|_{L^\infty}}2
\int_s^t
\|(w_1(\tau))_+\|_{L^\infty}\,d\tau.
\]

Setting
\[
M(t)=
\max\left\{
\|(w_1(t))_+\|_{L^\infty},
\|(w_2(t))_+\|_{L^\infty}
\right\},
\]
we obtain
\[
M(t)
\leq
\frac{\|d\|_{L^\infty}}2
\int_s^t M(\tau)\,d\tau.
\]
Hence Gronwall's inequality yields $M(t)=0$,
which proves the first assertion.

To prove \eqref{eq:comparison-maximum}, let
\[
M_0=
\max\left\{
\|w_1(s)\|_{L^\infty},
\|w_2(s)\|_{L^\infty}
\right\}.
\]
Since the constant functions $(M_0,M_0)$
and $(-M_0,-M_0)$ satisfy
\eqref{eq:comparison-system} with equality,
we can apply the first assertion to
$(w_1-M_0,w_2-M_0)$ and
$(-w_1-M_0,-w_2-M_0)$.
This gives
\[
-M_0\leq w_1(t,x),w_2(t,x)\leq M_0,
\]
and completes the proof.
\end{proof}

We next introduce the evolution operator
$\mathcal U_a(t,s)$ by using the characteristic
integral equations.

For given initial data $(f,g)\in C_b(\mathbb R)^2$
and $0\leq s\leq T$, we consider
\begin{align*}
y(t,x)
={}&
e^{-\frac12\int_s^t a(x_-(\xi))\,d\xi}
f(x_-(s))
\\
&+\frac12\int_s^t
e^{-\frac12\int_\tau^t a(x_-(\xi))\,d\xi}
a(x_-(\tau))z(\tau,x_-(\tau))\,d\tau,
\\
z(t,x)
={}&
e^{-\frac12\int_s^t a(x_+(\xi))\,d\xi}
g(x_+(s))
\\
&+\frac12\int_s^t
e^{-\frac12\int_\tau^t a(x_+(\xi))\,d\xi}
a(x_+(\tau))y(\tau,x_+(\tau))\,d\tau,
\end{align*}
where $x_\pm(\tau)=x_\pm(\tau;t,x)$
are the characteristic curves associated with $c(t,x)$.

These integral equations admit a unique bounded
continuous solution $(y,z)$ by successive approximations.
For $0\leq s\leq t\leq T$, we define
the evolution operator $\mathcal U_a(t,s)$ by
\[
\mathcal U_a(t,s)
\begin{pmatrix}
f\\g
\end{pmatrix}
=
\begin{pmatrix}
y(t,\cdot)\\z(t,\cdot)
\end{pmatrix}.
\]

If $(y,z)$ is sufficiently smooth,
the integral equations are equivalent to
the following linear transport system:
\begin{equation}
\begin{cases}
\displaystyle
\partial_-y=\frac{a(x)}2(z-y),
\\[2mm]
\displaystyle
\partial_+z=\frac{a(x)}2(y-z).
\end{cases}
\label{eq:homogeneous-transport}
\end{equation}

Writing $W=(y,z)^T$, this system can also
be expressed as
\[
\partial_tW=\mathcal L_a(t)W,
\]
where
\[
\mathcal L_a(t)
=
\begin{pmatrix}
c(t,x)\partial_x-\dfrac{a(x)}2
&
\dfrac{a(x)}2
\\[2mm]
\dfrac{a(x)}2
&
-c(t,x)\partial_x-\dfrac{a(x)}2
\end{pmatrix}.
\]

The evolution operator $\mathcal U_a(t,s)$
is a bounded linear operator on
$C_b(\mathbb R)^2$.
It satisfies
\[
\mathcal U_a(s,s)=I
\]
and
\[
\mathcal U_a(t,\tau)\mathcal U_a(\tau,s)
=
\mathcal U_a(t,s)
\]
for $s\leq\tau\leq t$.

By Lemma \ref{lem:comparison},
the evolution operator preserves nonnegative functions.
Moreover, comparison with constant solutions gives
\begin{equation}
\|\mathcal U_a(t,s)F\|_\infty
\leq\|F\|_\infty.
\label{eq:evolution-contraction}
\end{equation}

In what follows, transport equations with continuous solutions are understood in the sense of their characteristic integral equations, and all integrals involving the evolution operator are taken pointwise in $x$.
\begin{lemma}[Uniform estimate for the linear transport system]
\label{lem:linear-estimate}

Suppose that
\[
a(x)\geq0,
\qquad
a\in C_b^1(\mathbb R)\cap L^1(\mathbb R),
\]
and that \eqref{eq:linear-speed} holds.

Then, for every bounded continuous vector-valued
function $F$ on $[s,T]\times\mathbb R$,
\begin{equation}
\left\|
\int_s^t
\mathcal U_a(t,\tau)
\bigl[a(\cdot)F(\tau,\cdot)\bigr]
\,d\tau
\right\|_\infty
\leq
C_a
\sup_{s\leq\tau\leq t}
\|F(\tau)\|_\infty,
\label{eq:linear-uniform-estimate}
\end{equation}
where
\[
C_a
=
2\left(
\exp\left(\frac{\|a\|_{L^1}}{c_*}\right)-1
\right).
\]

In particular, the constant $C_a$ is independent
of $s$, $t$, and $T$.

\end{lemma}

\begin{proof}
To apply Lemma \ref{lem:comparison}, we construct a bounded nonnegative supersolution that is independent of time.

By linearity, we may assume that
\[
\sup_{s\leq\tau\leq t}\|F(\tau)\|_\infty\leq1.
\]
Writing $F=(f,g)$, we consider
\begin{equation}
\begin{cases}
\displaystyle
\partial_-y=\frac{a(x)}2(z-y)+a(x)f,\\[2mm]
\displaystyle
\partial_+z=\frac{a(x)}2(y-z)+a(x)g,\\[2mm]
y(s,x)=z(s,x)=0.
\end{cases}
\label{eq:forced-transport}
\end{equation}

We seek bounded supersolutions $P(x),Q(x)$
independent of time. Since $|f|,|g|\leq1$,
it is sufficient to find nonnegative functions
$P,Q$ satisfying
\[
\begin{cases}
\displaystyle
-c(t,x)P'(x)+\frac{a(x)}2(P-Q)\geq a(x),\\[2mm]
\displaystyle
c(t,x)Q'(x)+\frac{a(x)}2(Q-P)\geq a(x).
\end{cases}
\]

We look for such functions in the form
\[
P(x)=H+h(x),\qquad Q(x)=H-h(x),
\]
where $H$ is a nonnegative constant and $h\geq0$.
The above inequalities then become
\[
-ch'+ah\geq a,\qquad
-ch'-ah\geq a.
\]
Since $c(t,x)\geq c_*$, both inequalities
hold if
\[
h'(x)=-\frac{a(x)}{c_*}(1+h(x)).
\]
Solving this equation with $h(x)\to0$
as $x\to\infty$, we obtain
\[
h(x)=
\exp\left(
\frac1{c_*}\int_x^\infty a(\xi)\,d\xi
\right)-1.
\]

Since $a\geq0$ and $a\in L^1(\mathbb R)$,
\[
0\leq h(x)\leq H,
\qquad
H=\exp\left(\frac{\|a\|_{L^1}}{c_*}\right)-1.
\]
Thus $P,Q\geq0$ and
\[
\|P\|_{L^\infty},\|Q\|_{L^\infty}\leq2H.
\]
Moreover, the construction gives
\[
\partial_-P-\frac a2(Q-P)\geq a,
\qquad
\partial_+Q-\frac a2(P-Q)\geq a.
\]

Hence $(P,Q)$ is a supersolution of
\eqref{eq:forced-transport}.
Applying Lemma \ref{lem:comparison}, we obtain
\[
y(t,x)\leq P(x),\qquad z(t,x)\leq Q(x).
\]
Similarly, $(-P,-Q)$ is a subsolution, and thus
\[
-P(x)\leq y(t,x)\leq P(x),
\qquad
-Q(x)\leq z(t,x)\leq Q(x).
\]
Consequently,
\[
\|(y(t),z(t))\|_\infty
\leq2H
=
2\left(
e^{\|a\|_{L^1}/c_*}-1
\right).
\]

Since the solution with zero initial data
is represented by
\[
\begin{pmatrix}y(t)\\z(t)\end{pmatrix}
=
\int_s^t
\mathcal U_a(t,\tau)
\begin{pmatrix}
a(\cdot)f(\tau,\cdot)\\
a(\cdot)g(\tau,\cdot)
\end{pmatrix}
\,d\tau,
\]
the desired estimate follows.
\end{proof}

\subsection{Uniform estimates for the Riemann invariants}

We apply the preceding comparison argument to
\eqref{eq:riemann-system}.

\begin{lemma}
\label{lem:riemann-uniform}

Let $(u,v)$ be a classical solution of
\eqref{eq:Euler} with $u>0$.
Suppose that
\[
b(t)\geq0,\qquad a(x)\geq0.
\]

Then
\begin{equation}
\max\left\{
\|r(t)\|_{L^\infty},
\|s(t)\|_{L^\infty}
\right\}
\leq
\max\left\{
\|r(0)\|_{L^\infty},
\|s(0)\|_{L^\infty}
\right\}.
\label{eq:riemann-uniform}
\end{equation}

In particular, for sufficiently small
$\varepsilon>0$, we have
\begin{equation}
|u(t,x)-1|+|v(t,x)|
\leq C\varepsilon
\label{eq:uv-uniform}
\end{equation}
as long as the classical solution exists.

\end{lemma}

\begin{proof}

We set
\[
w_1=r,\qquad w_2=-s.
\]

Then \eqref{eq:riemann-system} becomes
\[
\begin{cases}
\displaystyle
\partial_-w_1
=\frac{b(t)+a(x)}2(w_2-w_1),
\\[2mm]
\displaystyle
\partial_+w_2
=\frac{b(t)+a(x)}2(w_1-w_2).
\end{cases}
\]

From Lemma \ref{lem:comparison},
comparison with constant solutions yields
\eqref{eq:riemann-uniform}.

From the definitions of $r$ and $s$,
we have
\[
v=\frac{r+s}{2}
\]
and
\[
\eta(u)
=
\frac{2}{\gamma-1}+\frac{s-r}{2}.
\]

Since
\[
\|r(0)\|_{L^\infty}
+\|s(0)\|_{L^\infty}
\leq C\varepsilon,
\]
we obtain
\[
|v(t,x)|\leq C\varepsilon
\]
and
\[
\left|
\eta(u(t,x))-\frac{2}{\gamma-1}
\right|
\leq C\varepsilon.
\]

The function $\eta$ is smooth and strictly
decreasing in a neighborhood of $u=1$.
Therefore, for sufficiently small $\varepsilon$,
\[
|u(t,x)-1|\leq C\varepsilon.
\]

This proves \eqref{eq:uv-uniform}.

\end{proof}

In particular, by taking $\varepsilon>0$
sufficiently small, we may use the uniform bound
\[
\frac12\leq c(u(t,x))\leq2
\]
throughout the existence interval.

\subsection{Transformations of the Riemann-invariant derivatives}

We introduce transformations of the derivatives
of the Riemann invariants which will be used
in the proof of Theorem \ref{thm:global}.

We denote by $c'(u)=dc(u)/du$ the derivative
of the sound speed with respect to $u$.

Differentiating \eqref{eq:riemann-system}
with respect to $x$, we obtain
\begin{equation}
\begin{aligned}
\partial_-r_x
={}&
\frac{c'(u)}{2c(u)}r_x(r_x-s_x)
-\frac{a+b}{2}(r_x+s_x)
-a'(x)v,
\\
\partial_+s_x
={}&
\frac{c'(u)}{2c(u)}s_x(s_x-r_x)
-\frac{a+b}{2}(r_x+s_x)
-a'(x)v.
\end{aligned}
\label{eq:riemann-derivative-system}
\end{equation}

\begin{remark}
Although the solution is assumed to be only $C^1$,
the derivatives $r_x$ and $s_x$ are differentiable
along their respective characteristic curves.
Indeed, this follows by differentiating the
characteristic integral equations for $r$ and $s$
with respect to the initial points of the characteristics.
This differentiation is justified since the characteristic
flows are $C^1$ with respect to their initial points.
Consequently, the above equations hold along
the corresponding characteristic curves.
\end{remark}

Here and below, we abbreviate
$c=c(u(t,x))$, $c'=c'(u(t,x))$,
$a=a(x)$, and $b=b(t)$
when no confusion arises.

To eliminate the terms involving $a'(x)v$,
we define
\[
R=r_x-\frac{a(x)}{c(u)}v,
\qquad
S=s_x+\frac{a(x)}{c(u)}v.
\]

Indeed, from \eqref{eq:Euler} and
\eqref{eq:riemann-derivatives}, we have
\[
\begin{aligned}
\partial_-v
&=v_t-cv_x
=-cs_x-(a+b)v,\\
\partial_+v
&=v_t+cv_x
=cr_x-(a+b)v.
\end{aligned}
\]
Since $\partial_-u=s_x$ and $\partial_+u=r_x$,
we obtain
\begin{equation}
\begin{aligned}
\partial_-\left(\frac acv\right)
={}&
-a'v-as_x
-\frac{a(a+b)}c v
-\frac{ac'}{c^2}vs_x,
\\
\partial_+\left(\frac acv\right)
={}&
a'v+ar_x
-\frac{a(a+b)}c v
-\frac{ac'}{c^2}vr_x.
\end{aligned}
\label{eq:correction-derivatives}
\end{equation}

Thus, substituting
\eqref{eq:correction-derivatives} into
\eqref{eq:riemann-derivative-system},
we obtain
\begin{equation}
\begin{aligned}
\partial_-R
={}&
\frac a2(S-R)-\frac b2(R+S)
+\frac{ab}{c}v
+\mathcal N_R,
\\
\partial_+S
={}&
\frac a2(R-S)-\frac b2(R+S)
-\frac{ab}{c}v
+\mathcal N_S,
\end{aligned}
\label{eq:RS-system}
\end{equation}
where
\[
\mathcal N_R
=
\frac{c'}{2c}R(R-S)
+\frac{ac'}{2c^2}v(3R+S)
\]
and
\[
\mathcal N_S
=
-\frac{c'}{2c}S(R-S)
-\frac{ac'}{2c^2}v(R+3S).
\]

We note that the derivatives of $a(x)$
no longer appear in \eqref{eq:RS-system}.

Next, to eliminate the cross terms involving
$b(t)$ in \eqref{eq:RS-system}, we introduce
\begin{equation}
\widetilde R
=
R+\frac{b(t)}2(u-1),
\qquad
\widetilde S
=
S+\frac{b(t)}2(u-1).
\label{eq:modified-RS}
\end{equation}

Since
\[
\partial_-(u-1)
=
s_x=S-\frac acv
\]
and
\[
\partial_+(u-1)
=
r_x=R+\frac acv,
\]
we obtain
\begin{equation}
\begin{cases}
\displaystyle
\partial_-\widetilde R
=
\frac a2(\widetilde S-\widetilde R)
-\frac b2\widetilde R+F_R,
\\[2mm]
\displaystyle
\partial_+\widetilde S
=
\frac a2(\widetilde R-\widetilde S)
-\frac b2\widetilde S+F_S,
\end{cases}
\label{eq:modified-RS-system}
\end{equation}
where
\[
F_R
=
\left(\frac{b'}2+\frac{b^2}{4}\right)(u-1)
+\frac{ab}{2c}v+\mathcal N_R
\]
and
\[
F_S
=
\left(\frac{b'}2+\frac{b^2}{4}\right)(u-1)
-\frac{ab}{2c}v+\mathcal N_S.
\]

Under the uniform estimate
\eqref{eq:uv-uniform}, we have
\[
|\mathcal N_R|+|\mathcal N_S|
\leq
C(|R|+|S|)^2
+C\varepsilon a(x)(|R|+|S|).
\]

The system \eqref{eq:modified-RS-system}
will be used to derive uniform estimates
for the derivatives of the solution.

\subsection{Time-dependent weights and linear estimates}

We finally establish estimates involving
the time-dependent damping coefficient
$b(t)$.

Recall that
\[
b(t)=\lambda(1+t)^{-\mu},
\]
where either
\[
\lambda>0,\quad 0<\mu<1,
\]
or
\[
\lambda>2,\quad \mu=1.
\]

For $0\leq s\leq t$, we define
\[
K_b(t,s)
=
\exp\left(
-\frac12\int_s^t b(\tau)\,d\tau
\right).
\]

\begin{lemma}
\label{lem:time-weight}

Under the assumptions on $\lambda$ and $\mu$,
there exists a constant $C>0$ such that
\begin{equation}
K_b(t,0)\leq Cb(t),
\label{eq:Kb-initial}
\end{equation}
\begin{equation}
\int_0^t K_b(t,s)b(s)^2\,ds
\leq Cb(t),
\label{eq:Kb-convolution}
\end{equation}
and
\begin{equation}
K_b(t,s)b(s)\leq Cb(t)
\label{eq:Kb-ratio}
\end{equation}
for all $0\leq s\leq t$.

Furthermore,
\begin{equation}
|b'(t)|\leq Cb(t)^2.
\label{eq:b-derivative}
\end{equation}

\end{lemma}

\begin{proof}

First, we have
\[
b'(t)
=
-\mu\lambda(1+t)^{-\mu-1}.
\]
Thus
\[
\frac{|b'(t)|}{b(t)^2}
=
\frac{\mu}{\lambda}(1+t)^{\mu-1}
\leq\frac{\mu}{\lambda},
\]
which proves \eqref{eq:b-derivative}.

We first consider the case $0<\mu<1$.
In this case,
\[
K_b(t,s)
=
\exp\left[
-\frac{\lambda}{2(1-\mu)}
\left\{
(1+t)^{1-\mu}-(1+s)^{1-\mu}
\right\}
\right].
\]

In particular, the stretched exponential decay
implies
\[
K_b(t,0)\leq C(1+t)^{-\mu}
\leq Cb(t).
\]

Next, we set
\[
J(t)=\int_0^tK_b(t,s)b(s)^2\,ds.
\]
Then
\[
J'(t)+\frac{b(t)}2J(t)=b(t)^2,
\qquad J(0)=0.
\]

Since
\[
\frac{|b'(t)|}{b(t)^2}
=
\frac{\mu}{\lambda}(1+t)^{\mu-1}
\longrightarrow0
\qquad(t\to\infty),
\]
there exists $T_0>0$ such that
\[
b'(t)+\frac12b(t)^2
\geq\frac14b(t)^2
\]
for $t\geq T_0$.

Consequently, if $C$ is sufficiently large,
the function $Cb(t)$ satisfies
\[
(Cb(t))'
+\frac{b(t)}2Cb(t)
\geq b(t)^2
\]
for $t\geq T_0$, and
\[
J(T_0)\leq Cb(T_0).
\]

The comparison principle for scalar
ordinary differential equations yields
\[
J(t)\leq Cb(t)
\qquad(t\geq T_0).
\]
Increasing $C$ if necessary, the estimate
also holds for $0\leq t\leq T_0$.
This proves \eqref{eq:Kb-convolution}.

We next show \eqref{eq:Kb-ratio}.
Since
\[
\frac{K_b(t,s)b(s)}{b(t)}
=
\exp\left(
\int_s^t
\left\{
\frac{\mu}{1+\tau}-\frac{b(\tau)}2
\right\}d\tau
\right),
\]
and
\[
\frac{\mu}{1+\tau}-\frac{b(\tau)}2
\leq0
\]
for all sufficiently large $\tau$,
the ratio is uniformly bounded.
Thus \eqref{eq:Kb-ratio} follows.

We next consider the case $\mu=1$ and
$\lambda>2$.
Then
\[
K_b(t,s)
=
\left(\frac{1+s}{1+t}\right)^{\lambda/2}.
\]

Since $\lambda/2>1$,
\[
K_b(t,0)
=(1+t)^{-\lambda/2}
\leq(1+t)^{-1}
=\frac{b(t)}{\lambda}.
\]

Moreover,
\begin{align*}
\int_0^tK_b(t,s)b(s)^2\,ds
&=
\frac{\lambda^2}{(1+t)^{\lambda/2}}
\int_0^t(1+s)^{\lambda/2-2}\,ds
\\
&=
\frac{2\lambda}{\lambda-2}b(t)
\left(
1-(1+t)^{1-\lambda/2}
\right)
\\
&\leq
\frac{2\lambda}{\lambda-2}b(t).
\end{align*}

Finally,
\[
\frac{K_b(t,s)b(s)}{b(t)}
=
\left(
\frac{1+s}{1+t}
\right)^{\lambda/2-1}
\leq1.
\]

This completes the proof.

\end{proof}

We now combine the preceding estimates with
the evolution operator $\mathcal U_a(t,s)$.

Consider the inhomogeneous system
\begin{equation}
\begin{cases}
\displaystyle
\partial_-y
=
\frac{a(x)}2(z-y)-\frac{b(t)}2y+f,
\\[2mm]
\displaystyle
\partial_+z
=
\frac{a(x)}2(y-z)-\frac{b(t)}2z+g.
\end{cases}
\label{eq:damped-linear-system}
\end{equation}

Writing $W=(y,z)$ and $F=(f,g)$,
we have
\[
\partial_tW
=
\left(
\mathcal L_a(t)-\frac{b(t)}2I
\right)W+F.
\]

Since $b(t)$ is independent of $x$,
its integrating factor is scalar.
Therefore, Duhamel's formula gives
\begin{equation}
\begin{aligned}
W(t)
={}&
K_b(t,s)\mathcal U_a(t,s)W(s)
\\
&+
\int_s^t
K_b(t,\tau)
\mathcal U_a(t,\tau)F(\tau)\,d\tau.
\end{aligned}
\label{eq:damped-duhamel}
\end{equation}

This representation remains valid when
$c=c(t,x)$ depends on both $t$ and $x$.

\begin{lemma}[Weighted estimate for the linear transport system]
\label{lem:weighted-linear}

Suppose that the assumptions of
Lemma \ref{lem:linear-estimate} and
Lemma \ref{lem:time-weight} are satisfied.

Then, for every bounded continuous vector-valued
function $G$ on $[0,T]\times\mathbb R$,
\[
\begin{aligned}
&\left\|
\int_0^t
K_b(t,s)\mathcal U_a(t,s)
\bigl[a(\cdot)b(s)G(s,\cdot)\bigr]\,ds
\right\|_\infty
\\
&\qquad\leq
Cb(t)
\sup_{0\leq s\leq t}\|G(s)\|_\infty,
\end{aligned}
\]
where $C$ is independent of $t$ and $T$.

\end{lemma}

\begin{proof}

We fix $t\in[0,T]$.

By Lemma \ref{lem:time-weight},
\[
\sup_{0\leq s\leq t}
\frac{K_b(t,s)b(s)}{b(t)}
\leq C.
\]

We define
\[
G_t(s,x)
=
\frac{K_b(t,s)b(s)}{b(t)}G(s,x),
\qquad 0\leq s\leq t.
\]

Then
\[
\sup_{0\leq s\leq t}\|G_t(s)\|_\infty
\leq
C\sup_{0\leq s\leq t}\|G(s)\|_\infty.
\]

Since $\mathcal U_a(t,s)$ is linear and
$K_b(t,s)b(s)/b(t)$ is independent of $x$,
we have
\begin{align*}
&\int_0^t
K_b(t,s)\mathcal U_a(t,s)
\bigl[a(\cdot)b(s)G(s)\bigr]\,ds
\\
&\qquad=
b(t)\int_0^t
\mathcal U_a(t,s)
\bigl[a(\cdot)G_t(s)\bigr]\,ds.
\end{align*}

Applying Lemma \ref{lem:linear-estimate}
to $G_t$, we obtain
\begin{align*}
&\left\|
\int_0^t
K_b(t,s)\mathcal U_a(t,s)
\bigl[a(\cdot)b(s)G(s)\bigr]\,ds
\right\|_\infty
\\
&\qquad\leq
C_ab(t)
\sup_{0\leq s\leq t}\|G_t(s)\|_\infty
\\
&\qquad\leq
Cb(t)
\sup_{0\leq s\leq t}\|G(s)\|_\infty.
\end{align*}

This completes the proof.

\end{proof}

\begin{corollary}
\label{cor:linear-source}

Consider the linear system
\eqref{eq:damped-linear-system} with initial data
$W(0)=W_0$.
Suppose that its forcing term can be written as
\[
F(t,x)
=
b(t)^2F_1(t,x)
+
a(x)b(t)F_2(t,x),
\]
where $F_1$ and $F_2$ are bounded and continuous on $[0,T]\times\mathbb R$.

Then
\begin{equation}
\begin{aligned}
\|W(t)\|_\infty
\leq Cb(t)\biggl(
&\|W_0\|_\infty+
\sup_{0\leq s\leq t}\|F_1(s)\|_\infty
\\
&+
\sup_{0\leq s\leq t}\|F_2(s)\|_\infty
\biggr)
\end{aligned}
\label{eq:linear-source-estimate}
\end{equation}
for all $t\in[0,T]$.

The constant $C$ is independent of $t$ and $T$.

\end{corollary}

\begin{proof}

By Duhamel's formula
\eqref{eq:damped-duhamel},
\begin{align*}
W(t)
={}&
K_b(t,0)\mathcal U_a(t,0)W_0
\\
&+
\int_0^t
K_b(t,s)\mathcal U_a(t,s)
[b(s)^2F_1(s)]\,ds
\\
&+
\int_0^t
K_b(t,s)\mathcal U_a(t,s)
[a(\cdot)b(s)F_2(s)]\,ds.
\end{align*}

The first term is estimated by
\eqref{eq:evolution-contraction} and
\eqref{eq:Kb-initial}.
The second term is estimated by
\eqref{eq:evolution-contraction} and
\eqref{eq:Kb-convolution}.
The last term is estimated by
Lemma \ref{lem:weighted-linear}.
Combining these estimates, we obtain
\eqref{eq:linear-source-estimate}.

\end{proof}

\section{Proof of the main theorem}

In this section, we prove Theorem \ref{thm:global}
by using the estimates established in Section 2.

\subsection{Uniform estimates}

Let $(u,v)$ be the unique local classical solution
of \eqref{eq:Euler}, and let $T^*>0$ denote its
maximal existence time.

We first establish uniform estimates for $u$ and $v$.
From the assumptions on the initial data,
we have
\[
\|r(0)\|_{L^\infty}
+
\|s(0)\|_{L^\infty}
\leq C\varepsilon.
\]

Since
\[
b(t)\geq0,\qquad a(x)\geq0,
\]
Lemma \ref{lem:riemann-uniform} implies that
\begin{equation}
\|r(t)\|_{L^\infty}
+
\|s(t)\|_{L^\infty}
\leq C\varepsilon
\label{eq:proof-rs-uniform}
\end{equation}
for all $t\in[0,T^*)$.

Consequently,
\begin{equation}
\|u(t)-1\|_{L^\infty}
+
\|v(t)\|_{L^\infty}
\leq C\varepsilon.
\label{eq:proof-uv-uniform}
\end{equation}

In particular, if $\varepsilon>0$ is sufficiently
small, then
\begin{equation}
\frac12\leq c(u(t,x))\leq2
\label{eq:proof-speed}
\end{equation}
for all $(t,x)\in[0,T^*)\times\mathbb R$.

Moreover, since $u$ remains in a compact
neighborhood of $1$, we have
\begin{equation}
|c(u)|+|c(u)^{-1}|+|c'(u)|
\leq C.
\label{eq:proof-c-bounds}
\end{equation}

The constants in
\eqref{eq:proof-rs-uniform}--\eqref{eq:proof-c-bounds}
are independent of $t$ and $T^*$.

\subsection{Estimates for $\widetilde R$ and $\widetilde S$}

We next establish uniform estimates for the
modified variables introduced in
\eqref{eq:modified-RS}.
Recall that
\[
R=r_x-\frac{a(x)}{c(u)}v,
\qquad
S=s_x+\frac{a(x)}{c(u)}v
\]
and
\[
\widetilde R
=
R+\frac{b(t)}2(u-1),
\qquad
\widetilde S
=
S+\frac{b(t)}2(u-1).
\]

We set
\[
W(t,x)
=
\begin{pmatrix}
\widetilde R(t,x)\\
\widetilde S(t,x)
\end{pmatrix}.
\]

By \eqref{eq:modified-RS-system}, $W$ satisfies the following system in the sense of the characteristic integral equations:
$W$ satisfies
\begin{equation}
\partial_tW
=
\left(
\mathcal L_a(t)-\frac{b(t)}2I
\right)W+F,
\label{eq:proof-W-system}
\end{equation}
where
\[
F=
\begin{pmatrix}
F_R\\
F_S
\end{pmatrix}
\]
and
\begin{align*}
F_R
={}&
\left(\frac{b'}2+\frac{b^2}{4}\right)(u-1)
+\frac{ab}{2c}v
+\mathcal N_R,
\\
F_S
={}&
\left(\frac{b'}2+\frac{b^2}{4}\right)(u-1)
-\frac{ab}{2c}v
+\mathcal N_S.
\end{align*}

Here
\begin{align}
\mathcal N_R
={}&
\frac{c'}{2c}R(R-S)
+\frac{ac'}{2c^2}v(3R+S),
\label{eq:proof-NR}
\\
\mathcal N_S
={}&
-\frac{c'}{2c}S(R-S)
-\frac{ac'}{2c^2}v(R+3S).
\label{eq:proof-NS}
\end{align}

We shall prove that
\[
\|W(t)\|_\infty
\leq C\varepsilon b(t)
\]
for all $t\in[0,T^*)$.

Fix $T\in(0,T^*)$ arbitrarily, and set
\[
M(T)
=
\sup_{0\leq t\leq T}
\frac{\|W(t)\|_\infty}{\varepsilon b(t)}.
\]

Since $b(t)>0$ and the solution is classical,
$M(T)$ is finite for every $T<T^*$.

From the definitions of $\widetilde R$ and
$\widetilde S$, we have
\[
R=\widetilde R-\frac b2(u-1),
\qquad
S=\widetilde S-\frac b2(u-1).
\]

Therefore, using \eqref{eq:proof-uv-uniform},
we obtain
\begin{equation}
|R(t,x)|+|S(t,x)|
\leq
C\varepsilon b(t)(1+M(T))
\label{eq:proof-RS-bound}
\end{equation}
for all $(t,x)\in[0,T]\times\mathbb R$.

We now estimate the nonlinear terms.

By \eqref{eq:proof-NR},
\eqref{eq:proof-NS},
\eqref{eq:proof-c-bounds}, and
\eqref{eq:proof-RS-bound}, we have
\begin{align}
|\mathcal N_R|+|\mathcal N_S|
\leq{}&
C\varepsilon^2b(t)^2(1+M(T))^2
\nonumber\\
&+
C\varepsilon^2a(x)b(t)(1+M(T)).
\label{eq:proof-nonlinear}
\end{align}

On the other hand, Lemma \ref{lem:time-weight}
gives
\[
|b'(t)|\leq Cb(t)^2.
\]

Hence, from \eqref{eq:proof-uv-uniform},
we obtain
\[
\left|
\left(\frac{b'}2+\frac{b^2}{4}\right)(u-1)
\right|
\leq C\varepsilon b(t)^2.
\]

Similarly,
\begin{equation}
\left|\frac{a(x)b(t)}{2c(u)}v\right|
\leq C\varepsilon a(x)b(t).
\label{eq:proof-space-source}
\end{equation}

Combining
\eqref{eq:proof-nonlinear}--\eqref{eq:proof-space-source},
we obtain
\begin{align*}
|F_R|+|F_S|
\leq{}&
C\varepsilon b(t)^2
+
C\varepsilon a(x)b(t)
\\
&+
C\varepsilon^2b(t)^2(1+M(T))^2
\\
&+
C\varepsilon^2a(x)b(t)(1+M(T)).
\end{align*}

We now write the forcing term in the form
\begin{equation}
F(t,x)
=
b(t)^2F_1(t,x)
+
a(x)b(t)F_2(t,x).
\label{eq:proof-source-decomposition}
\end{equation}

More precisely, we define
\[
F_1=
\begin{pmatrix}
F_{1,R}\\
F_{1,S}
\end{pmatrix},
\qquad
F_2=
\begin{pmatrix}
F_{2,R}\\
F_{2,S}
\end{pmatrix},
\]
where
\begin{align*}
F_{1,R}
={}&
\left(
\frac{b'}{2b^2}+\frac14
\right)(u-1)
+
\frac{c'}{2cb^2}R(R-S),
\\
F_{1,S}
={}&
\left(
\frac{b'}{2b^2}+\frac14
\right)(u-1)
-
\frac{c'}{2cb^2}S(R-S),
\end{align*}
and
\begin{align*}
F_{2,R}
={}&
\frac{v}{2c}
+
\frac{c'}{2c^2b}v(3R+S),
\\
F_{2,S}
={}&
-\frac{v}{2c}
-
\frac{c'}{2c^2b}v(R+3S).
\end{align*}

From the preceding estimates, it follows that
\begin{equation}
\sup_{0\leq t\leq T}\|F_1(t)\|_\infty
\leq
C\varepsilon
+
C\varepsilon^2(1+M(T))^2
\label{eq:proof-F1-bound}
\end{equation}
and
\begin{equation}
\sup_{0\leq t\leq T}\|F_2(t)\|_\infty
\leq
C\varepsilon
+
C\varepsilon^2(1+M(T)).
\label{eq:proof-F2-bound}
\end{equation}

The constants in these estimates are independent
of $T$.

We next estimate the initial value of $W$.

Since
\[
u(0,x)=1+\varepsilon\varphi(x),
\qquad
v(0,x)=\varepsilon\psi(x),
\]
we have
\[
u_t(0,x)=v_x(0,x)
=\varepsilon\psi_x(x).
\]

Moreover,
\[
c(u(0,x))u_x(0,x)
=
\varepsilon c(u(0,x))\varphi_x(x).
\]

It follows from the definitions of $R,S$ that
\[
\|R(0)\|_{L^\infty}
+
\|S(0)\|_{L^\infty}
\leq C\varepsilon.
\]

Since $b(0)=\lambda$, we obtain
\begin{equation}
\|W(0)\|_\infty\leq C\varepsilon.
\label{eq:proof-W-initial}
\end{equation}

We now apply Corollary \ref{cor:linear-source}
to \eqref{eq:proof-W-system}.

The assumptions of that corollary are satisfied
because \eqref{eq:proof-speed} gives a uniform
positive lower bound for the characteristic speed,
and the forcing term has the decomposition
\eqref{eq:proof-source-decomposition}.

Consequently, for $0\leq t\leq T$,
\begin{align*}
\|W(t)\|_\infty
\leq{}&
Cb(t)\biggl(
\|W(0)\|_\infty
+
\sup_{0\leq s\leq t}\|F_1(s)\|_\infty
\\
&\hspace{32mm}
+
\sup_{0\leq s\leq t}\|F_2(s)\|_\infty
\biggr).
\end{align*}

Using
\eqref{eq:proof-F1-bound},
\eqref{eq:proof-F2-bound}, and
\eqref{eq:proof-W-initial}, we find
\begin{align*}
\|W(t)\|_\infty
\leq{}&
C_0\varepsilon b(t)
\\
&+
C_1\varepsilon^2b(t)(1+M(T))^2
\\
&+
C_2\varepsilon^2b(t)(1+M(T)).
\end{align*}

Dividing by $\varepsilon b(t)$,
taking the supremum over $0\leq t\leq T$,
and increasing the constants if necessary,
we obtain
\begin{equation}
M(T)
\leq
C_0+C_1\varepsilon(1+M(T))^2.
\label{eq:proof-bootstrap-inequality}
\end{equation}

We now close this estimate by a continuity argument.

Choose a constant $K>0$ such that
\[
K>4(C_0+1)
\]
and
\[
M(0)<\frac K2.
\]

We take $\varepsilon>0$ sufficiently small
so that
\[
C_1\varepsilon(1+K)^2<\frac K4.
\]

Suppose that $M(T)\leq K$.
Then \eqref{eq:proof-bootstrap-inequality} gives
\[
M(T)
\leq
C_0+C_1\varepsilon(1+K)^2
<
\frac K2.
\]

Thus the bootstrap bound improves strictly.

Since the solution and its characteristic
derivatives depend continuously on time on
every compact subinterval of $[0,T^*)$,
the usual continuity argument yields
\[
M(T)\leq K
\qquad\text{for every }T<T^*.
\]

The constants $K$ and $\varepsilon$ can be
chosen independently of $T$.

It follows from the definition of $M(T)$ that
\[
\|\widetilde R(t)\|_{L^\infty}
+
\|\widetilde S(t)\|_{L^\infty}
\leq C\varepsilon b(t)
\]
for all $t\in[0,T^*)$.

Finally, since
\[
R=\widetilde R-\frac b2(u-1),
\qquad
S=\widetilde S-\frac b2(u-1),
\]
we conclude from \eqref{eq:proof-uv-uniform} that
\begin{equation}
\|R(t)\|_{L^\infty}
+
\|S(t)\|_{L^\infty}
\leq C\varepsilon b(t).
\label{eq:proof-RS-decay}
\end{equation}

\subsection{Conclusion of the proof}

We now derive estimates for the derivatives of
the original unknown functions.

Recall that
\[
R=r_x-\frac{a(x)}{c(u)}v,
\qquad
S=s_x+\frac{a(x)}{c(u)}v.
\]

From \eqref{eq:riemann-derivatives}, we have
\[
u_t=v_x=\frac{R+S}{2}.
\]

Using \eqref{eq:proof-RS-decay},
we obtain
\[
\|u_t(t)\|_{L^\infty}
+
\|v_x(t)\|_{L^\infty}
\leq C\varepsilon b(t).
\]

On the other hand, we have
\[
u_x
=
\frac{r_x-s_x}{2c(u)}
=
\frac{R-S}{2c(u)}
+\frac{a(x)}{c(u)^2}v.
\]

Consequently,
\[
|u_x(t,x)|
\leq
C\varepsilon\bigl(b(t)+a(x)\bigr).
\]

Furthermore, from the second equation of
\eqref{eq:Euler}, we have
\[
\begin{aligned}
v_t
&=c(u)^2u_x-(a(x)+b(t))v\\
&=\frac{c(u)}2(R-S)-b(t)v.
\end{aligned}
\]

Hence, by \eqref{eq:proof-uv-uniform} and
\eqref{eq:proof-RS-decay},
we obtain
\[
\|v_t(t)\|_{L^\infty}
\leq C\varepsilon b(t).
\]

Recalling that
\[
b(t)=\lambda(1+t)^{-\mu},
\]
we conclude that
\begin{equation}
\|u_t(t)\|_{L^\infty}
+
\|v_t(t)\|_{L^\infty}
+
\|v_x(t)\|_{L^\infty}
\leq
C\varepsilon(1+t)^{-\mu}
\label{eq:proof-final-time}
\end{equation}
and
\begin{equation}
|u_x(t,x)|
\leq
C\varepsilon
\left(
(1+t)^{-\mu}+a(x)
\right).
\label{eq:proof-final-space}
\end{equation}

Since $a\in C_b^1(\mathbb R)$,
\eqref{eq:proof-final-space} also implies
\[
\sup_{0\leq t<T^*}
\|u_x(t)\|_{L^\infty}<\infty.
\]

Combining these estimates with
\eqref{eq:proof-uv-uniform}, we obtain
\[
\begin{aligned}
\sup_{0\leq t<T^*}
\bigl(
&\|u(t)\|_{L^\infty}
+\|v(t)\|_{L^\infty}
\\
&+\|u_t(t)\|_{L^\infty}
+\|v_t(t)\|_{L^\infty}
\\
&+\|u_x(t)\|_{L^\infty}
+\|v_x(t)\|_{L^\infty}
\bigr)<\infty.
\end{aligned}
\]

Moreover, \eqref{eq:proof-uv-uniform}
guarantees that $u(t,x)$ remains uniformly
bounded away from zero.
Thus the continuation criterion for the
quasilinear hyperbolic system
\eqref{eq:Euler} implies that the maximal
existence time cannot be finite.
Therefore,
\[
T^*=\infty.
\]
The uniqueness of the global solution follows
from the local uniqueness theorem.
The estimates
\eqref{eq:zeroth-estimate},
\eqref{eq:time-derivative-estimate}, and
\eqref{eq:space-derivative-estimate}
follow from
\eqref{eq:proof-uv-uniform},
\eqref{eq:proof-final-time}, and
\eqref{eq:proof-final-space}, respectively.
This completes the proof of
Theorem \ref{thm:global}.

\section*{Declarations}

\subsection*{Funding}

The author was partially supported by Grants-in-Aid for Scientific Research
(C) (No.~23K03169) and Grant-in-Aid for Challenging Research (Pioneering)
(No.~24K21198).

\subsection*{Competing interests}

The author has no competing interests to declare that are relevant to the
content of this article.

\subsection*{Data availability}

Data sharing is not applicable to this article as no datasets were generated
or analyzed during the current study.

\subsection*{Declaration of AI use}
 During the preparation of this manuscript, the author used ChatGPT (OpenAI) to assist with mathematical discussions and the revision of the manuscript. All mathematical results and proofs were independently verified by the author, who takes full responsibility for the content of the manuscript.

\bibliographystyle{plain}
\bibliography{biblio}

@article{Friedrichs1948,
  author  = {Kurt Otto Friedrichs},
  title   = {Nonlinear hyperbolic differential equations
             for functions of two independent variables},
  journal = {American Journal of Mathematics},
  volume  = {70},
  number  = {3},
  pages   = {555--589},
  year    = {1948}
}

@article{Lax1964,
  author  = {Peter David Lax},
  title   = {Development of Singularities of Solutions of
             Nonlinear Hyperbolic Partial Differential Equations},
  journal = {Journal of Mathematical Physics},
  volume  = {5},
  number  = {5},
  pages   = {611--613},
  year    = {1964},
  doi     = {10.1063/1.1704154}
}

@article{HsiaoLiu1992,
  author  = {Ling Hsiao and Tai-Ping Liu},
  title   = {Convergence to nonlinear diffusion waves
             for solutions of a system of hyperbolic
             conservation laws with damping},
  journal = {Communications in Mathematical Physics},
  volume  = {143},
  number  = {3},
  pages   = {599--605},
  year    = {1992},
  doi     = {10.1007/BF02099268}
}

@article{HsiaoLiu1993,
  author  = {Ling Hsiao and Tai-Ping Liu},
  title   = {Nonlinear diffusive phenomena of
             nonlinear hyperbolic systems},
  journal = {Chinese Annals of Mathematics, Series B},
  volume  = {14},
  number  = {4},
  pages   = {465--480},
  year    = {1993}
}

@article{HsiaoSerre1996,
  author  = {Ling Hsiao and Denis Serre},
  title   = {Global Existence of Solutions for the System
             of Compressible Adiabatic Flow through Porous Media},
  journal = {SIAM Journal on Mathematical Analysis},
  volume  = {27},
  number  = {1},
  pages   = {70--77},
  year    = {1996},
  doi     = {10.1137/S0036141094267078}
}

@article{Nishihara1997,
  author  = {Kenji Nishihara},
  title   = {Asymptotic Behavior of Solutions of Quasilinear
             Hyperbolic Equations with Linear Damping},
  journal = {Journal of Differential Equations},
  volume  = {137},
  number  = {2},
  pages   = {384--395},
  year    = {1997},
  doi     = {10.1006/jdeq.1997.3268}
}

@article{MarcatiNishihara2003,
  author  = {Pierangelo Marcati and Kenji Nishihara},
  title   = {The {$L^p$}--{$L^q$} estimates of solutions
             to one-dimensional damped wave equations
             and their application to the compressible
             flow through porous media},
  journal = {Journal of Differential Equations},
  volume  = {191},
  number  = {2},
  pages   = {445--469},
  year    = {2003},
  doi     = {10.1016/S0022-0396(03)00026-3}
}

@article{Mei2010,
  author  = {Ming Mei},
  title   = {Best Asymptotic Profile for Hyperbolic
             {$p$}-System with Damping},
  journal = {SIAM Journal on Mathematical Analysis},
  volume  = {42},
  number  = {1},
  pages   = {1--23},
  year    = {2010},
  doi     = {10.1137/090756594}
}

@article{Wirth2004,
  author  = {Jens Wirth},
  title   = {Solution representations for a wave equation
             with weak dissipation},
  journal = {Mathematical Methods in the Applied Sciences},
  volume  = {27},
  number  = {1},
  pages   = {101--124},
  year    = {2004},
  doi     = {10.1002/mma.446}
}

@article{Wirth2006,
  author  = {Jens Wirth},
  title   = {Wave equations with time-dependent dissipation.
             {I}. Non-effective dissipation},
  journal = {Journal of Differential Equations},
  volume  = {222},
  number  = {2},
  pages   = {487--514},
  year    = {2006},
  doi     = {10.1016/j.jde.2005.07.019}
}

@article{Wirth2007,
  author  = {Jens Wirth},
  title   = {Wave equations with time-dependent dissipation.
             {II}. Effective dissipation},
  journal = {Journal of Differential Equations},
  volume  = {232},
  number  = {1},
  pages   = {74--103},
  year    = {2007},
  doi     = {10.1016/j.jde.2006.06.004}
}

@article{Pan2016Global,
  author  = {Xinghong Pan},
  title   = {Global existence of solutions to {1-D}
             {Euler} equations with time-dependent damping},
  journal = {Nonlinear Analysis},
  volume  = {132},
  pages   = {327--336},
  year    = {2016},
  doi     = {10.1016/j.na.2015.11.022}
}

@article{Pan2016Blowup,
  author  = {Xinghong Pan},
  title   = {Blow up of solutions to {1-D}
             {Euler} equations with time-dependent damping},
  journal = {Journal of Mathematical Analysis and Applications},
  volume  = {442},
  number  = {2},
  pages   = {435--445},
  year    = {2016},
  doi     = {10.1016/j.jmaa.2016.04.075}
}

@article{ChenPanZhu2017,
  author  = {Geng Chen and Ronghua Pan and Shengguo Zhu},
  title   = {Singularity Formation for the
             Compressible {Euler} Equations},
  journal = {SIAM Journal on Mathematical Analysis},
  volume  = {49},
  number  = {4},
  pages   = {2591--2614},
  year    = {2017},
  doi     = {10.1137/16M1062818}
}

@article{ChenYoungZhang2013,
  author  = {Geng Chen and Robin Young and Qingtian Zhang},
  title   = {Shock formation in the compressible {Euler}
             equations and related systems},
  journal = {Journal of Hyperbolic Differential Equations},
  volume  = {10},
  number  = {1},
  pages   = {149--172},
  year    = {2013},
  doi     = {10.1142/S0219891613500069}
}

@article{Sugiyama2018,
  author  = {Yuusuke Sugiyama},
  title   = {Singularity formation for the {1D}
             compressible {Euler} equations with
             variable damping coefficient},
  journal = {Nonlinear Analysis},
  volume  = {170},
  pages   = {70--87},
  year    = {2018},
  doi     = {10.1016/j.na.2017.12.013}
}

@incollection{Sugiyama2020,
  author    = {Yuusuke Sugiyama},
  title     = {Remark on global existence of solutions
               to the {1D} compressible {Euler} equation
               with time-dependent damping},
  booktitle = {Advanced Studies in Pure Mathematics},
  volume    = {85},
  pages     = {379--389},
  year      = {2020},
  publisher = {Mathematical Society of Japan},
  doi       = {10.2969/aspm/08510379}
}

@article{ChenLiLiMeiZhang2020,
  author  = {Shaohua Chen and Haitong Li and Jingyu Li
             and Ming Mei and Kaijun Zhang},
  title   = {Global and blow-up solutions for compressible
             {Euler} equations with time-dependent damping},
  journal = {Journal of Differential Equations},
  volume  = {268},
  number  = {9},
  pages   = {5035--5077},
  year    = {2020},
  doi     = {10.1016/j.jde.2019.11.002}
}

@article{GengLinMei2020,
  author  = {Shifeng Geng and Yanping Lin and Ming Mei},
  title   = {Asymptotic Behavior of Solutions to {Euler}
             Equations with Time-Dependent Damping
             in Critical Case},
  journal = {SIAM Journal on Mathematical Analysis},
  volume  = {52},
  number  = {2},
  pages   = {1463--1488},
  year    = {2020},
  doi     = {10.1137/19M1272846}
}

@article{Chen2024,
  author  = {Shaohua Chen},
  title   = {Global solutions for the {1-D} compressible
             {Euler} equations with time-dependent damping},
  journal = {Communications on Pure and Applied Analysis},
  volume  = {23},
  number  = {1},
  pages   = {1--20},
  year    = {2024},
  doi     = {10.3934/cpaa.2023125}
}

@article{GengLaiYuenZhou2023,
  author  = {Jinbo Geng and Ning-An Lai
             and Man Wai Yuen and Jiang Zhou},
  title   = {Blow-up for compressible {Euler} system
             with space-dependent damping in {1-D}},
  journal = {Advances in Nonlinear Analysis},
  volume  = {12},
  number  = {1},
  note    = {Article 20220304},
  year    = {2023},
  doi     = {10.1515/anona-2022-0304}
}

@article{Sugiyama2025,
  author  = {Yuusuke Sugiyama},
  title   = {Global existence and blow-up for the
             {1D} damped compressible {Euler} equations
             with time and space dependent perturbation},
  journal = {Nonlinear Analysis},
  volume  = {250},
  note    = {Article 113658},
  year    = {2025},
  doi     = {10.1016/j.na.2024.113658}
}

@article{CaiLiuMeiWang2025,
  author  = {Yang Cai and Changchun Liu
             and Ming Mei and Zejia Wang},
  title   = {Asymptotic behavior of solutions to the
             {Cauchy} problem for {1D} {$p$}-system
             with spatiotemporal damping:
             Case 1. {$v_+=v_-$}},
  journal = {Journal of Differential Equations},
  volume  = {436},
  note    = {Article 113347},
  year    = {2025},
  doi     = {10.1016/j.jde.2025.113347}
}

@article{CaiLiuMeiWang2026,
  author  = {Yang Cai and Changchun Liu
             and Ming Mei and Zejia Wang},
  title   = {Novel convergence of solutions to {1D}
             compressible {Euler} equations with
             spatiotemporal damping in critical case},
  journal = {Journal of Differential Equations},
  volume  = {461},
  note    = {Article 114094},
  year    = {2026},
  doi     = {10.1016/j.jde.2026.114094}
}

@article{Pan2021,
  author  = {Xinghong Pan},
  title   = {Global existence and asymptotic behavior
             of solutions to the Euler equations
             with time-dependent damping},
  journal = {Applicable Analysis},
  volume  = {100},
  number  = {16},
  pages   = {3546--3575},
  year    = {2021},
  doi     = {10.1080/00036811.2020.1722805}
}

\end{document}